\documentclass[11pt]{amsart}
\usepackage{latexsym,amssymb,amsmath}
\usepackage{graphics}
\usepackage{mathptmx} 
\begin{document}
\sloppy
\renewcommand{\theequation}{\arabic{section}.\arabic{equation}}
\thinmuskip = 0.5\thinmuskip
\medmuskip = 0.5\medmuskip
\thickmuskip = 0.5\thickmuskip
\arraycolsep = 0.3\arraycolsep

\newtheorem{theorem}{Theorem}[section]
\newtheorem{lemma}[theorem]{Lemma}
\newtheorem{prop}[theorem]{Proposition}
\newtheorem{definition}{Definition}
\newtheorem{remark}{Remark}
\renewcommand{\thetheorem}{\arabic{section}.\arabic{theorem}}
\newcommand{\prf}{\noindent{\bf Proof.}\ }
\def\prfe{\hspace*{\fill} $\Box$

\smallskip \noindent}

\def\be{\begin{equation}}
\def\ee{\end{equation}}
\def\bea{\begin{eqnarray}}
\def\eea{\end{eqnarray}}
\def\beas{\begin{eqnarray*}}
\def\eeas{\end{eqnarray*}}
\def\x{{\bf x}}
\def\S{\mathcal{S}}

\newcommand{\R}{\mathbb R}
\newcommand{\N}{\mathbb N}
\newcommand{\T}{\mathbb T}
\newcommand{\K}{\mathbb S}
\newcommand{\E}{\mathcal E}
\newcommand{\D}{\mathcal D}
\newcommand{\U}{\mathcal U}
\newcommand{\q}{\mathbb P}
\newcommand{\z}{{\bf P}}
\def\a{{\bf a}}
\def\g{\partial}

\title{Orthogonality conditions and asymptotic stability in the Stefan problem
        with surface tension}
\author{Mahir Had\v{z}i\'{c}} 
%
%
%
\maketitle
\begin{abstract}
We prove nonlinear asymptotic stability of steady spheres in the two-phase Stefan problem with surface tension.
Our method relies on the introduction of appropriate orthogonality conditions in conjunction
with a high-order energy method.
\end{abstract}
\section{Introduction}\label{se:intro}
We are interested in the question of long-time nonlinear stability of steady state solutions 
to the two-phase Stefan problem with surface tension, one of the best known 
parabolic free boundary problems. 
It is a simple model of phase transitions in liquid-solid systems.
\par
Let $\Omega\subset\R^n$ denote a $C^1$-domain that contains a liquid and a solid separated by an interface
$\Gamma$. As the melting or cooling take place the boundary moves and we are naturally led
to a free boundary problem. 
Define the solid phase $\Omega^-(t)$ as a region encircled by $\Gamma(t)$ and define the
liquid phase $\Omega^+(t):=\Omega\setminus{\overline{\Omega^-}}$.
The unknowns are the location of the interface $\{\Gamma(t);\,t\geq0\}$ and the temperature function
$
v:[0,T]\times\Omega\to\R.
$
Let $\Gamma_0$ be the initial position 
of the free boundary and $v_0:\Omega\to\R$ be the initial temperature.
We denote the normal velocity of
$\Gamma$ by $V$ and normalize it to be {\it positive if $\Gamma$ is locally expanding $\Omega^+(t)$}.
Furthermore, we denote the mean curvature of $\Gamma$ by $\kappa$.
With these notations, $(v,\Gamma)$ satisfies the following free boundary value problem:
\be\label{eq:temp}
\partial_t v-\Delta v=0\quad\textrm{in}\quad\Omega\setminus\Gamma.
\ee
\be\label{eq:dirichlet}
v=\kappa\quad\text{on}\quad\Gamma,
\ee
\be\label{eq:neumann}
V=[v_n]^+_-\quad\textrm{on}\quad\Gamma,
\ee
\be\label{eq:neumann1}
v_n=0\quad\textrm{on}\quad\g\Omega,
\ee
\be\label{eq:initial}
v(0,\cdot)=v_0;\quad\Gamma(0)=\Gamma_0.
\ee
Given $v$, we write
$v^+$ and $v^-$ for the restriction of $v$ to $\Omega^+(t)$ and $\Omega^-(t)$,
respectively. With this notation
$[v_n]^+_-$ stands for the jump
of the normal derivatives across the interface $\Gamma(t)$, namely
$
[v_n]^+_-:=v_n^+-v^-_n,
$
where $n$ stands for the unit normal on the hypersurface $\Gamma(t)$ with respect to $\Omega^+(t)$.
In~(\ref{eq:neumann1}) $\g\Omega$ stands for the outer fixed boundary of $\Omega$.
Two basic identities related to the above problem are the ``mass" conservation law:
\be\label{eq:mass}
\g_t\Big[\int_{\Omega}v(t,x)\,dx+|\Omega^-(t)|\Big]=0,
\ee
and the energy dissipation law:
\[
\g_t\Big[\frac{1}{2}\int_{\Omega}v^2\,dx+|\Gamma(t)|\Big]+\int_{\Omega}|\nabla v|^2\,dx=0.
\]
Here $|\Omega^-(t)|$ and $|\Gamma(t)|$ stand for the Lebesgue volume of 
$\Omega^-(t)$ and the surface area
of $\Gamma(t)$ respectively.
Steady states of the above problem consist of static spheres and they form an 
$(n+1)$-dimensional family $\mathcal{F}$:
\be\label{eq:steady}
\mathcal{F}:=\big\{\Sigma(R,\a)\big|\,\,\,\,
\a\in\R^n,\,R\in\R_+\big\},
\ee
where for any $\a=(a^1,\dots,a^n)\in\R^n$, $R\in\R_+$ the pair
\[
\Sigma(R,\a):=((n-1)/R, S_R(\a))
\] 
is a time-independent solution of the problem~(\ref{eq:temp}) -~(\ref{eq:neumann1}),
if $S_R(\a)\subset\Omega$. 
Let us parametrize the moving boundary $\Gamma$ as a graph over a given steady state 
$\Sigma(R,\a)\in\mathcal{F}$: we introduce the radius function 
$r:[0,\infty[\times\K^{n-1}\to\R$ and the parametrization
$\phi:[0,\infty[\times\K^{n-1}\to\Gamma$ such that
\be\label{eq:par}
\phi(t,\xi)=\a+r(t,\xi)\xi,\quad\xi\in\K^{n-1}
\ee
and define the perturbation $(u,\,f):=(v-(n-1)/R,\,r-1)$.
To $(u,f)$ we associate a high-order energy norm given by
\[
||(u,f)(t)||=\sum_{k=0}^N\big\{||u||_{W^{N-k,\infty}W^{2k,2}([0,t]\times\Omega)}^2
+||\nabla_gf||_{W^{N-k,\infty}W^{2k+1,2}([0,t]\times\K^{n-1})}^2\big\}.
\]
Our main result is the following:
\begin{theorem}\label{th:main}
Assume that $\zeta_R>0$, where
$
\zeta_R:=\frac{1}{|\Omega|}-\frac{n-1}{|S_R|R^2}
$
(i.e. the linear stability criterion holds).
Then there exists $\epsilon>0$, such that if initially 
\[
||(u,f)(0)||<\epsilon,
\]
where
$(u,f)(t):=(v(t,\cdot)-(n-1)/R,r(t,\cdot)-R)$
is the perturbation from the steady state $\Sigma(R,\a)$,
then there exists a global-in-time unique solution $(v,\Gamma)(t)$ to the 
problem~(\ref{eq:temp}) -~(\ref{eq:initial}). Moreover, there exists a pair $(\bar{R},\bar{\a})$
close to $(R,\a)$ such that $(v,\Gamma)$ converges exponentially fast to 
$\Sigma(\bar{R},\bar{\a})\in\mathcal{F}$. More precisely, if we parametrize
$(v,\Gamma)(t)$ as a perturbation of $\Sigma(\bar{R},\bar{\a})$, by setting
$\Gamma(t)=\big\{\bar{\a}+\bar{r}(t,\xi)\xi,\,\,\,\,\xi\in\K^{n-1}\big\}$ and define
$(\bar{u},\bar{f})(t)=(v(t,\cdot)-(n-1)/\bar{R},\,\bar{r}(t,\cdot)-\bar{R})$ then
there exist constants $C_1$ and $C_2$ such that
\[
||(\bar{u},\bar{f})(t)||\leq C_1e^{-C_2t},\quad t\geq0.
\]
\end{theorem}
\par
{\bf Remark.}
By choosing $N$ large enough in the definition of the norm $||\cdot||$, 
we can make the solution as smooth as desired in a classical sense.
We have not aimed for finding the minimal $N$; rather, for the sake of clarity and
conciseness of the estimates, we allow ourselves the flexibility of keeping $N$ sufficiently large. 
\par
{\bf Remark.}
Note that the time derivatives occurring in the expression $||(u,f)(0)||$ are implicitly given
by terms with only spatial derivatives, via the equations~(\ref{eq:temp}) and~(\ref{eq:neumann}).
\subsection{Notation}
\label{su:notation}
On the unit sphere $\K^{n-1}$ the Riemannian gradient with respect to the standard metric is
denoted by $\nabla_g$ and the Laplace-Beltrami operator by $\Delta_g$. For a given function $h$
the $k$-th time derivative is interchangeably denoted by $\g_{t^k}h$ or $h_{t^k}$. When writing
various norms of the 
functions defined on the unit sphere $\K^{n-1}$, we drop the domain from the notation, e.g.
$||f||_{L^2}:=||f||_{L^2(\K^{n-1})}$. A ball (sphere) of radius $R$ centered at a point $\a\in\R^n$
will be denoted by $B_R(\a)$ ($S_R(\a)$). An important tool in our analysis are the spherical harmonics
on the unit sphere.
For a given function $h\in L^2(\K^{n-1})$,
we introduce the spherical harmonics decomposition for $h$
\[
h=\sum_{k=0}^{\infty}\sum_{i=1}^{m(k)}h_{k,i}s_{k,i}.
\] 
Here for each $k\in\N_0$, $s_{k,i}$, $i=1,\dots,m(k)$ stand for the spherical harmonics of degree $k$.
They are defined as restrictions of homogeneous polynomials of degree $k$ in $n$ variables onto
the unit sphere. The set 
$\bigcup_{k=0}^{\infty}\bigcup_{i=1}^{m(k)}\big\{s_{k,i}\big\}$
forms an orthonormal basis on $\K^{n-1}$ with respect to $L^2$-product,
thus justifying the above expansion. For more details, we refer the reader to~\cite{Mu}.
The first $n+1$ spherical harmonics $s_{0,1},s_{1,1},\dots,s_{1,n}$ will be denoted by
$s_0,s_1,\dots,s_n$. For $i=1,\dots,n$ we define the first momenta of $h$ $\q_1h_i$ 
as well as the mean $\q_0h$ by setting
\be
\q_1h_i:=\int_{\K^{n-1}}hs_i\,d\xi;\quad\q_0h=\frac{1}{|\K^{n-1}|}\int_{\K^{n-1}}h\,d\xi.
\ee
We denote $\q_{2+}h:=h-\q_0h-\sum_{k=1}^n\q_1h_i$. For a given function $\U$ we interchangeably denote
$\g_{x^i}\U=\U_{x^i}=\U_i$ and similarly for mixed partial derivatives. Einstein's summation convention
is employed and we sum over repeated indices. A generic constant appearing in various estimates is
denoted by $C$ and it may change from line to line. The unit normal to a given surface 
$\Gamma$ is denoted by $n_{\Gamma}$ and the index is dropped if the surface in question is clear from the 
context.
\subsection{Previous work}
Stefan problem has been studied in a variety of mathematical literature over the past century
(see for instance \cite{Vi}).
If~(\ref{eq:dirichlet}) is replaced by the boundary condition $v=0$ on $\Gamma$, 
the resulting problem is called the {\it classical} 
Stefan problem.
It has been known that the classical Stefan problem admits unique global weak solutions in several dimensions
(\cite{Fri1}, \cite{Fri2} and \cite{Ka}). 
The references to the regularity of weak solutions of the two-phase classical Stefan problem
are, among others, ~\cite{AtCaSa1}, ~\cite{AtCaSa2}, ~\cite{CaEv}, ~\cite{Di}. 
Local classical solutions are established in~\cite{Han} and~\cite{Me}. 
\par
If the diffusion equation~(\ref{eq:temp}) is
replaced by the elliptic equation $\Delta u=0$, then
the resulting problem is called the Mullins-Sekerka problem
(also known as the quasi-stationary Stefan problem or the Hele-Shaw problem with surface tension).  
Global existence for
the two-phase quasi-stationary Stefan problem close to a sphere
in two dimensions has been obtained in~\cite{Ch} and~\cite{CoPu}, and in arbitrary dimensions in~\cite{EsSi1}.
Global stability for the one-phase quasi-stationary Stefan problem is established
in~\cite{FriRe2}. Local-in-time
solutions in parabolic H\"older spaces in arbitrary dimensions are established in~\cite{ChHoYi}.
\par
Global stability of a single equilibrium for the Mullins-Sekerka problem follows 
also from the work of \textsc{Alikakos \& Fusco}~\cite{AlFu},
as an immediate consequence of their study of Ostwald ripening. See the remarks at the end of Section~\ref{se:intro} and
Section~\ref{se:he} respectively. 
\par
In a recent preprint~\cite{PrSiZa} the authors analyze the thermodynamically consistent Stefan problem with surface tension.
In the stable regime $\zeta_R>0$ their methods imply the stability of stationary spheres for the Stefan problem with surface tension using
the maximal regularity techniques in $L^p$-type parabolic Sobolev spaces.
\par
As to the Stefan problem with surface tension (also known as the Stefan problem with Gibbs-Thomson correction),
global weak existence theory (without uniqueness) is developed in~\cite{AlWa},~\cite{Lu} and~\cite{Roe}.
In \cite{FriRe} the authors consider the
Stefan problem with small surface tension i.e. $\sigma\ll1$ whereby~(\ref{eq:dirichlet})
is replaced by $v=\sigma\kappa$.
Local existence for the Stefan problem is studied in \cite{Ra}.
In \cite{EsPrSi} the authors
prove a local existence and uniqueness result 
under a smallness assumption on the initial datum
close to flat hypersurfaces. Linear stability and instability results for spheres 
are contained in~\cite{PrSi}. 

The first global-in-time nonlinear stability result for the 
flat steady hypersurfaces was given in~\cite{HaGu1}. 
Some of the references for the Stefan problem with surface tension and 
kinetic undercooling effects are~\cite{ChRe},~\cite{Ra},~\cite{So},~\cite{Vi2}. 
\subsection{Motivation, methods and plan of the paper}
To explain the linear stability criterion $\zeta_R>0$ and motivate our result, let us look at the linearization around a fixed steady state 
$\Sigma(1,\a)\in\mathcal{F}$:
\be\label{eq:linear}
\g_t(u,f)=\mathcal{L}(u,f),
\ee
where $\mathcal{L}(u,f)=(\Delta u,-[\g_nu]^+_-)$, $u=-(n-1)f-\Delta_gf$ on $\K^{n-1}$.
This linearization is easily obtained if we parametrize the
moving boundary $\Gamma$ as a graph over $S_1(\a)$: with the
radius function $r=1+f$ and the parametrization $\phi$ as in~(\ref{eq:par}), the
perturbation $(u,f)$ takes the form $(u,\,f)=(v-(n-1),\,R-1)$. It is readily checked that 
up to the first order $\kappa\circ\phi=(n-1)-(n-1)f-\Delta_gf$
and $V\circ\phi=f_t$. Associated to the linear problem~(\ref{eq:linear}) is
the energy dissipation identity (cf.~\cite{HaGu2}):
\[
\g_t\Big\{\int_{\Omega}u^2+\int_{\K^{n-1}}\{|\nabla_gf|^2-(n-1)f^2\}\Big\}=-2\int_{\Omega}|\nabla u|^2.
\]
Note that it is not even clear whether the above energy is positive definite, due to
the presence of the negative definite term $-\int_{\K^{n-1}}(n-1)f^2$. 
Indeed, problem~(\ref{eq:linear}) may allow for a strictly positive eigenvalue under certain 
assumptions on the relative size of the domain $\Omega$ and the steady state $\Sigma(R,\a)$.
More precisely, to any $\Sigma(R,\a)\in\mathcal{F}$ we associate the stability parameter $\zeta_R$:
\[
\zeta_R:=\frac{1}{|\Omega|}-\frac{n-1}{|S_R|R^2},
\]
where $|\Omega|$ and $|S_R|$ stand for the Lebesgue volume of $\Omega$ and the surface area
of $S_R$ respectively. It turns out that the sign of $\zeta_R$
is critical to the stability properties of $\Sigma(R,\a)$: 
if $\zeta_R>0$ the solution is {\em linearly} stable and if $\zeta_R<0$ it is 
{\em linearly} unstable (\cite{HaGu2},~\cite{PrSi}). 
The full  {\em nonlinear} instability in the case $\zeta_R<0$,
under the additional assumption that $\Omega$ is a perfect ball of a given radius $R_*$ is proved in~\cite{HaGu2}. 
A closer look at the {\em instability} proof in~\cite{HaGu2} shows that the proof 
itself is insensitive to the specific shape of $\Omega$ and works 
for general domains. On the other hand, in the stability regime ($\zeta_R>0$) the situation
is more complicated.
To explain this, note that the linearized problem~(\ref{eq:linear}) has 
$(n+1)$ non-decaying solutions $\{\sigma_i\}_{i=0,\dots,n}$, 
being exactly the $(n+1)$ eigenvectors spanning the null-space of $\mathcal{L}$:
\[
\sigma_0:=(n-1,-1),\,\,\,\sigma_i=(0,s_i)\,\,\,\,i=1,\dots, n,
\]
where $s_i,\,i=1,\dots,n$, are the spherical harmonics defined in Subsection~\ref{su:notation}.
The existence of $\sigma_i$-s, $i=0,\dots,n$, encodes the $(n+1)$-dimensionality 
of the set $\mathcal{F}$.  We expect the perturbation of $\Sigma(R,\a)$ to converge
to a nearby asymptotic state $\Sigma(\bar{R},\bar{\a})$. On the other hand,
note that the first modes $\q_1f_i=\int_{\K^{n-1}}fs_i\,d\xi s_i$ ($i=1,\dots,n$)
in the spherical-harmonics expansion of $f$ are precisely cancelled by
the expression $Z(h):=\int_{\K^{n-1}}\{|\nabla_gh|^2-(n-1)h^2\}$ and they cannot be a-priori controlled 
by the energy. 
This feature of the problem is another manifestation of the non-trivial null-space structure 
of the linearized operator $\mathcal{L}$.
It causes a major difficulty in controlling the first momenta $\q_1f_i$, $i=1,\dots,n$,
thus introducing additional analytic difficulties.

In this paper, as stated in Theorem~\ref{th:main}, we will prove the {\em nonlinear asymptotic
stability} of the steady state $\Sigma(R,\a)\in\mathcal{F}$ in the
stability regime $\zeta_R>0$ for {\em general} domains $\Omega$. Our method has 
two basic ingredients. We start 
by introducing a set of geometrically motivated orthogonality conditions, that allow us to 
``mod out" an inherent degeneracy related to
the existence of the non-trivial nullspace of the linearized operator.
In analytical terms, we express $(v,\Gamma)$
in a set of new, tubular coordinates.
The second ingredient
is the high-order energy method developed in~\cite{HaGu1},~\cite{HaGu2}
as a part of the program to investigate stability and instability of steady states
in phase transition phenomena. 
The goal is to prove an energy estimate of the form
\be\label{eq:main}
\E(t)+\int_0^t\D(\tau)\,d\tau\leq\E(0)+\big(\delta+
C\sup_{0\leq\tau\leq t}\sqrt{\E(\tau)}\big)\int_0^t\D(\tau)\,d\tau,
\ee
where $\E$ and $\D$ are the {\em energy} and {\em dissipation} naturally associated to the problem
(cf.~(\ref{eq:en}) and~(\ref{eq:di})).
If $\E$ and $\delta>0$ are small, we can absorb the right-most term above into LHS, to obtain 
an a-priori estimate
\[
\E(t)+\frac{1}{2}\int_0^t\D(\tau)\,d\tau\leq\E(0),
\]
thus recovering the smallness assumption on $\E$ if $\E(0)$ is small. Through an iteration scheme 
and a continuity argument, this reasoning is made rigorous: proof of~(\ref{eq:main})
is based on a series of energy estimates for the error terms. However, while in
the previous works~\cite{HaGu1} and~\cite{HaGu2} the boundary
conditions~(\ref{eq:dirichlet}) and~(\ref{eq:neumann}) honored the energy method fully,
by introducing the tubular coordinates new error terms of (only) quadratic order emerge
(Lemma~\ref{lm:energy}). 
The precise structure of the
linearized curvature operator is exploited to get cancellation for such terms and to finally close the estimates
(Theorem~\ref{th:local}). 
To obtain global existence for the solution $(v,\Gamma)$, we use a Poincar\'e-inequality type
bound $\E\leq C\D$. It allows us to prove some decay of the solution on a bounded time
interval. Together with a suitable smallness
assumption on the initial data this decay can be bootstrapped to yield a global-in-time existence result.
\par
We wish to point out that our method can also be used to prove asymptotic stability of steady spheres for
the Mullins-Sekerka problem in arbitrary dimensions.
\par
{\bf Remark.}
After the publication of this work, the author learned that the stability of Mullins-Sekerka steady sphere follows easily from the
study of Ostwald ripening~\cite{AlFu,AlFuKa} where the problem
of long-time behavior close to a steady state
consisting of a finite union of steady spheres is analyzed. 
Moreover, the same orthogonality conditions as in this paper are introduced in~\cite{AlFu,AlFuKa}
to deal with the non-trivial null-space - see the remark at the end of Section~\ref{se:he}.
\par
The plan of the paper is as follows. In Section~\ref{se:orth}
we heuristically motivate and then introduce the orthogonality conditions.
We then derive the modulation equations and reformulate
the Stefan problem in the set of new, tubular coordinates.
In Section~\ref{se:fix}, the Stefan problem is reformulated on a fixed domain and
the high-order energy is introduced. 
In Section~\ref{se:apriori} we prove a Poincar\'e type estimate (Lemma~\ref{lm:poincare})
and the positive-definiteness of the energy (Lemma~\ref{lm:positive}). With these preparations,
energy estimates are performed and local existence theorem is proved in Section~\ref{se:local}.
Finally, the proof of the main result (Theorem~\ref{th:main}) is presented in Section~\ref{se:global}.
\section{Orthogonality conditions. Evolution problem}\label{se:orth}
To facilitate the analysis, without loss of generality, we 
assume a few simplifications. We will perturb away from the 
steady state $\Sigma(1,{\bf 0})$ (assuming $S_1({\bf 0})\subset\Omega$).
Furthermore, note that the mass conservation law~(\ref{eq:mass}) necessarily determines the radius of the final
asymptotic state. Recalling the parametrization~(\ref{eq:par}) of $\Gamma(t)$,~(\ref{eq:mass})
takes the form $\g_tM(v,\Gamma)(t)=0$, where
\[
M(v,\Gamma)(t):=\int_{\Omega}v(t,x)\,dx
+\int_{\K^{n-1}}\frac{r(t,\xi)^n}{n}\,d\xi.
\]
We will henceforth assume that initially
\be\label{eq:initialmass}
M(v_0,\Gamma_0)=M(\Sigma(1,{\bf 0})).
\ee
Condition~(\ref{eq:initialmass}) forces the 
asymptotic steady state to have the radius $R_{asymp}=1$ due to the
conservation of $M(v,\Gamma)(t)$. This assumption constraints our stability 
analysis to the ``manifold" of steady states 
$\mathcal{G}\subset\mathcal{F}$ (defined in~(\ref{eq:steady})) consisting of the 
elements $\Sigma(\a):=\Sigma(1,\a)$ of the fixed radius $R=1$:
\[
\mathcal{G}:=\big\{\Sigma(\a)\big|\,\,\,\,\a\in\Omega,\,\,\,\,S_1(\a)\subset\Omega\big\}.
\]
Finally, we assume $\int_{\K^{n-1}}f_0s_i=0$, $i=1,\dots,n$ where the
initial surface $\Gamma_0$ is parametrized
by $(1+f_0(\xi))\xi$, $\xi\in\K^{n-1}$. Otherwise, translate the steady state to a nearby one, so that the 
condition is satisfied.
\subsection{Heuristics}\label{se:he}
To motivate the analysis of the present work, in the following we provide some geometric 
heuristics for the choice of the above mentioned orthogonality conditions. Namely,
let us think of $\mathcal{G}$ as a submanifold of the subspace $\mathcal{H}$ 
of $L^2(\Omega)\times L^2(\K^{n-1};\,\Omega)$ of the functions of the form $(u,\,\a+R(\xi)\xi)$:
\[
\mathcal{H}:=\big\{(u,\a+r(\xi)\xi),\,\,\,\a\in\R^n,\,\,u:\Omega\to\R,\,\,\,r:\K^{n-1}\to\R\big\}
\subset L^2(\Omega)\times L^2(\K^{n-1};\,\Omega).
\]
Thus, for $\a\in\R^n$, $\Sigma(\a)(x,\xi)=(n-1,\a+\xi)\in\mathcal{G}$. The non-decaying solutions 
$\{\sigma_i\}_{i=1,\dots,n}$ of~(\ref{eq:linear})
correspond to the infinitesimal changes in the center coordinate parameters $a_i$, $i=1,\dots,n$.
Motivated by this observation,  
we expect the solution to the full nonlinear problem to decompose into a component tangential to
$\mathcal{G}$ and the 
dissipating part that belongs to a plane {\it transversal} to $\mathcal{G}$ in $\mathcal{H}$.
In fact, we shall demand that this plane is exactly the fiber
$\mathcal{G}^{\perp}$, $L^2$-orthogonal to $\mathcal{G}$ in $\mathcal{H}$.  
In other words, we choose the tubular coordinates $(\a(t),u(t,\cdot),f(t,\cdot))$ 
in a neighborhood of the steady state manifold $\mathcal{G}$:
\be\label{eq:decomp}
(v(t),\Gamma(t))=\Sigma_{\a(t)}+(u(t),f(t,\xi)\xi),
\ee
such that 
\be\label{eq:orthgeom}
(u(t,\cdot),f(t,\cdot)\,\iota(\cdot))\in\mathcal{G}^{\perp},
\ee
where $\iota:\K^{n-1}\to\R^n$ is the inclusion operator.
Condition~(\ref{eq:orthgeom}) is 
precisely the {\it orthogonality} condition. To obtain an analytic expression
for~(\ref{eq:orthgeom}), we note that 
the tangent space $T_{\a(t)}\mathcal{G}$ is spanned by the set
$\{\Sigma_i(x,\xi)=\g_{a^i}\Sigma_{\a(t)}\big|_{a^i=a^i(t)}\}_{i=1,\dots,n}=\{(0,e_i)\}_{i=1,\dots,n}$.
Thus, condition~(\ref{eq:orthgeom}) implies that for any $i=1,\dots,n$:
\be\label{eq:ortho}
0=\big((u,f\,\iota),(0,e_i)\big)_{L^2}=\int_{\K^{n-1}}f(t,\xi)\xi\cdot e_i\,d\xi=\int_{\K^{n-1}}f(t,\xi)s_i\,d\xi.
\ee
Here $e_i$ stands for the $i$-th unit vector in $\R^n$.
We have provided a geometric picture that renders the right orthogonality 
condition~(\ref{eq:ortho}),
thus drawing a parallel to the work of {\sc Friesecke} \& {\sc Pego}~\cite{FrPe2} on the stability
of solitons in Fermi-Pasta-Ulam lattices. There, the authors introduce suitable orthogonality conditions
exploiting the symplectic structure of the linearized problem. Their perturbed solution belongs
to the fiber {\em symplectically} orthogonal to the soliton state manifold.
\par
{\bf Remark.}
After the publication of this work, the author learned that the same orthogonality conditions as in~(\ref{eq:ortho}) were {\em first} introduced and used by
\textsc{Alikakos} \& \textsc{Fusco}~\cite{AlFu} (see also~\cite{AlFuKa}) in the study of Ostwald ripening in the context of the related Mullins-Sekerka
problem (see in particular eqns. (0.24) and (0.25)
from~\cite{AlFu} and for the heuristic motivation see the remarks on top of the page 443 in~\cite{AlFu}).
Their methods rely on optimal regularity techniques and,
among other, imply the stability of a single equilibrium in the Mullins-Sekerka problem as pointed out in the introduction. 
Same idea is also used in a related work on the normalized mean curvature flow~\cite{AlFr}. 
\subsection{Modulation equations. The evolution problem}
A natural question raised by the form of the orthogonality condition~(\ref{eq:ortho}) is to 
find the modulation equation (i.e. evolution equation) for the vector $\a(t)$.
This is not straightforward 
because $\a(t)$ does not explicitly appear in~(\ref{eq:ortho}), thus preventing us from differentiating
in time directly. To resolve this difficulty, for 
any $i=1,\dots,n$, let us define the {\it momentum test functions} 
$p_i(t,\cdot)\in C^{\infty}(\Omega)$:
$
p_i(x,t)=x^i-a^i(t),
$
where $x=(x^1,\dots, x^n)$. 
Assume for the moment that $(v,\Gamma)(t)$ is a classical solution of the Stefan 
problem~(\ref{eq:temp}) -~(\ref{eq:initial}). We then multiply the equation~(\ref{eq:temp})
by the momentum test functions, use the integration by parts once and the boundary condition
$[\partial_n v]^+_-=V_{\Gamma}$. As a result, for $i=1,\dots, n$, we obtain the identity
\be\label{eq:orth}
\int_{\Gamma(t)}V_{\Gamma}p_i=\int_{\Omega}u_tp_i+\int_{\Omega}\nabla u\cdot\nabla p_i.
\ee
Note that the normal velocity expressed in local coordinates takes the form
\[
V_{\Gamma}\circ\phi(t,\xi)=-\frac{r_tr}{|g|}+
\dot{\a}\cdot \vec{n}_{\Gamma}\circ\phi,\qquad |g|=\sqrt{r^2+|\nabla_gr|^2}
\]
and the mean curvature $\kappa$ takes the form:
\[
\kappa\circ\phi(t,\xi)=\frac{n-1}{|g|}-\frac{1}{R}\nabla_g\cdot\frac{\nabla_gR}{|g|}
=(n-1)-(n-1)f-\Delta_gf+N(f),
\]
where $N(f)$ stands for the quadratic nonlinear remainder.
The volume element $dS(\Gamma)$ takes the form $r^{n-2}|g|\,d\xi$ in the local
coordinates on the sphere, so we obtain 
\beas
\int_{\Gamma(t)}V_{\Gamma}p_i\,dS(\Gamma)&=&
-\int_{\K^{n-1}}r^nr_t\xi^i+
\dot{\a}\cdot\int_{\Gamma(t)}p_i\vec{n}_{\Gamma}
=-\int_{\K^{n-1}}r^nr_t\xi^i+
\dot{\a}\cdot\int_{\Omega^-(t)}\nabla p_i\\
&=&-\int_{\K^{n-1}}r^nr_t\xi^i+\dot{a}^i(t)|\Omega^-(t)|,
\eeas
where we used the Stokes theorem in the second equality and $\nabla p_i=e_i$ in the last. 
Plugging this back into~(\ref{eq:orth}), we conclude
\be\label{eq:orthA}
|\Omega^-(t)|\dot{a}^i(t)=\int_{\K^{n-1}}r^nr_t\xi^i
+\int_{\Omega}u_tp_i+\int_{\Omega}u_{x^i}.
\ee
Thus, if $(v,\Gamma)$ is the solution of the problem with the above choice of
coordinate description of $\Gamma$, then a fortiori, the moving center components
$a^i(t)$ satisfy the differential equation~(\ref{eq:orthA}). 
Written in terms of $f$, for any $i=1,\dots,n$, the first term on RHS of~(\ref{eq:orthA}) takes the form:
\[
\int_{\K^{n-1}}r^nr_ts_i\,d\xi=\int_{\K^{n-1}}(1+f)^nf_ts_i\,d\xi
=\sum_{k=1}^n\int_{\K^{n-1}}{n\choose k}f^kf_ts_i,
\]
where in the second equality we used the orthogonality condition~(\ref{eq:orth}), 
implying in particular $\int_{\K^{n-1}}f_ts_i=0$, $i=1,\dots n$.
The unknowns to be solved for, are the perturbation $(u,f)$ and the moving
center $\a$. Setting ${\bf p}=(p_1,\dots,p_n)$,
Stefan problem with surface tension~(\ref{eq:temp}) -~(\ref{eq:initial}) 
in the tubular coordinates $(u(t,\cdot), f(t,\cdot), \a(t))$ takes the form
\be\label{eq:tempu}
\partial_t u-\Delta u=0\quad\textrm{in}\quad\Omega(t).
\ee
\be\label{eq:dirichletu}
u=\kappa-(n-1)\quad\textrm{on}\quad\Gamma(t),
\ee
\be\label{eq:neumannu}
[u_n]^+_-\circ\phi=
-\frac{r_tr}{|g|}+\dot{\a}\cdot\vec{n}_{\Gamma}\circ\phi.
\quad\textrm{on}\quad\K^{n-1},
\ee
\be\label{eq:neumann1u}
u_n=0\quad\textrm{on}\quad\g\Omega,
\ee
\be\label{eq:orthAB}
\dot{\a}(t)=\frac{1}{|\Omega^-(t)|}\sum_{k=1}^n\int_{\K^{n-1}}{n\choose k}f^kf_t\xi
+\frac{1}{|\Omega^-(t)|}\int_{\Omega(t)}u_t{\bf p}+\frac{1}{|\Omega^-(t)|}\int_{\Omega(t)}\nabla u.
\ee
\be\label{eq:conditionsu}
\int_{\K^{n-1}}f\xi={\bf 0}.
\ee
\be\label{eq:initialu}
u(0,\cdot)=u_0;\quad\Gamma(0)=\Gamma_0;
\quad \a(0)={\bf 0}.
\ee
\be\label{eq:initialmassu}
M(1+u_0,\Gamma_0)=M(\Sigma(1,{\bf 0})).
\ee
The idea is
to exploit dissipative properties of the 
problem~(\ref{eq:tempu}) -~(\ref{eq:neumann1u}) to obtain 
a time-decay estimate for $(u,f)$. Plugging that
into~(\ref{eq:orthAB}), we hope to obtain a decay estimate
for $\dot{\a}$. By bootstrapping this procedure, we will ``drive" the
solution to its asymptotic equilibrium.
\section{Fixing the domain and the energy}\label{se:fix}
We first describe the pull back of
the problem~(\ref{eq:tempu}) -~(\ref{eq:initialmassu}) onto the fixed domain
$\Omega\setminus\K^{n-1}$ and then define the high-order energies $\E$ and $\D$,
deriving in particular the corresponding energy identities. From this point onwards,
we will denote $\S:=\K^{n-1}$. Define the following change of variables:
\[
 \Theta(t,x)=\pi(t,x)(x-\a(t)),
\]
where $\pi$ is a smooth scalar-valued function with the following property:
\be
\label{eq:change1}
\pi(t,x)=\left\{
\begin{array}{l}
\frac{1}{r(\frac{\x-\a(t)}{|\x-\a(t)|})},\quad|\x-\a(t)|-1\leq d,\\
1,\quad|\x-\a(t)|-1\geq 2d
\end{array}
\right.
\quad
d\quad\mbox{is small}.
\ee
For any $\x\in\Gamma$, we may write $\x=\a(t)+r(\frac{\x-\a(t)}{|\x-\a(t)|})\frac{x-\a(t)}{|x-\a(t)|}=\a(t)+r(\xi)\xi$,
by the definition of $\Gamma$. Hence
\[
\Theta(t,\x)=\pi\big(t,\a(t)+r(\xi)\frac{\x-\a(t)}{|\x-\a(t)|}\big)r(\xi)\frac{\x-\a(t)}{|\x-\a(t)|}=
\frac{1}{r(\xi)}r(\xi)\frac{\x-\a(t)}{|\x-\a(t)|}
=\frac{\x-\a(t)}{|\x-\a(t)|}\in\S.
\]
Thus $\Theta$ does the job for us: $\Theta:\Omega\to\Omega\setminus\S$. Note 
that this map is natural from the point of view of our \emph{geometric} assumption
that the evolving interfaces $\Gamma$ are, in fact, 
graphs over the unit sphere $\S$.
\par
The inverse map $x(\bar{x})$ to the above change of variables~$\Theta$
is given by 
\[
x=\a(t)+\rho(t,\bar{x})\bar{x}.
\] 
To convince ourselves that the function $\rho(t,\cdot):\Omega\setminus\S\to\R$ is well defined, we observe
that $\rho$ has to satisfy the relation
\[
\pi(t,\a(t)+\rho(t,\bar{x})\bar{x})\rho(\bar{x})-1=0.
\]
The existence of such a $\rho$ can be established by the implicit function theorem, applied
to the equation $F(t,s,\bar{x})=0$, where $F:\R\times\R\times\R^n\to\R$ is given by
$
F(t,s,v)=\pi(t,\a+sv)s-1.
$
Namely, $F_s(t,s,v)=(\nabla\pi\cdot v)s+\pi>0$ since $\pi\approx1$ and 
$\nabla\pi=O(\nabla_gf)$ is small.
We define $w:\Omega\setminus\S\to\R$ by setting
\be\label{eq:uchange}
w(t,\bar{x}):=u(t,\a(t)+\rho(t,\bar{x})\bar{x}).
\ee
The heat operator $\partial_t-\Delta$ on the domain $\Omega$ will transform
into a more complicated operator in the new coordinates as stated in the following lemma:
\begin{lemma}\label{lm:push}
The push forward of the operator $\g_t-\Delta$ with respect to the map $\Theta:\Omega\to\S$
reads:
\[
(\g_t-\Delta)^{\#}u=w_t-a_{ij}w_{\bar{x}^i\bar{x}^j}-b_iw_{\bar{x}^i},
\]
where
\be\label{eq:ab}
a_{ij}:=\pi^2\delta_{ij}+2\pi x^j\pi_{x^i}+x^ix^j|\nabla\pi|^2,\quad
b_i:=2\pi_{x^i}+\Delta\pi x^i-\pi_tx^i.
\ee
Furthermore, 
\be\label{eq:neumanntheta}
[w_n]^+_-
=\frac{r^2}{|g|}\Big([u_n]^+_-\Big)\circ\phi
\ee
\end{lemma}
The proof of this lemma is a straightforward calculation and it is presented in the appendix.
In conclusion, on the fixed domain $\Omega\setminus\S$, 
equations~(\ref{eq:tempu}) -~(\ref{eq:neumannu}) take the form
\be\label{eq:tempw}
w_t-a_{ij}w_{ij}-b_iw_{i}=0\quad\mbox{in}\,\,\,\,\Omega,
\ee
\be
w=\kappa\circ\phi-(n-1)\quad\mbox{on}\,\,\,\,\S,
\ee
\be
[w_n]^+_-=\frac{r^2}{|g|}\big(-\frac{r_tr}{|g|}+\dot{\a}\cdot n_{\Gamma}\circ\phi\big)
\quad\mbox{on}\,\,\,\,\S.
\ee
The boundary condition~(\ref{eq:neumann1u}) reads $\g_nw=0$ on $\g\Omega$.
The modulation equation~(\ref{eq:orthAB}) is easily expressed in the fixed coordinates
\be
\dot{\a}(t)=\frac{1}{|\Omega^-(t)|}\Big(\sum_{k=1}^n\int_{\K^{n-1}}{n\choose k}f^kf_t\xi
+\int_{\Omega}(u_t{\bf p})\circ\Theta|\mbox{det}D\Theta^{-1}|
+\int_{\Omega}\nabla u\circ\Theta|\mbox{det}D\Theta^{-1}|\Big),
\ee
where formulas~(\ref{eq:timenew}) and~(\ref{eq:grad}) are used to express
$u_t\circ\Theta$ and $\nabla u\circ\Theta$ in terms of the function $w=u\circ\Theta$.
The orthogonality condition~(\ref{eq:conditionsu}) retains its form and
the initial conditions take the form:
\be\label{eq:conditionsfixed}
w(0,\cdot)=w_0:=u_0\circ\Theta;\quad f(0,\cdot)=f_0;\quad\a(0)={\bf 0}.
\ee
\subsection{The model problem}
Let $\mu:\Omega\to\R_+$
be a smooth non-negative cut-off function such that $\mu=0$ close to $\S$ and
$\mu=1$ close to $\partial\Omega$ and the origin
\footnote{By ``close" we mean in an open neighborhood of prescribed positive thickness.}. 
For any $i=1,\dots,n$ let us define a differential operator
\[
 D^i=\mu\partial_{x^i}+(1-\mu)\partial_{\xi^i},
\]
where $\partial_{\xi^i}$ is the partial tangential differentiation operator defined in Cartesian coordinates through
\[
\partial_{\xi^i}u=\partial_{x^i}u-\frac{x^i}{|x|}\nabla u\cdot\frac{x}{|x|}
\]
For given functions
$v,w\in H^1(\Omega)$ a simple integration by parts shows
\[
\int_{\Omega}D^ivw=-\int_{\Omega}vD^iw-\int_{\Omega}vw\nu(\mu)+\int_{\partial\Omega}vwn^i,
\]
where $\nu(\mu):=(\partial_{x^i}-\partial_{\xi^i})\mu=\frac{x^i}{|x|}\nabla\mu\cdot\frac{x}{|x|}$, thus giving us
the integration-by-parts formula for the operator $D^i$.
For a given multi-index $m=(m_1,\dots,m_n)$ and a non-negative natural number $s\in\N_0$ we define 
\[
 D^m_s:=\partial_t^s\g^{m}.
\]
If we denote 
$
 \mathcal{L}:=a_{ij}\partial_{x^ix^j}+b_i\partial_{x^i}
$
the second order elliptic operator appearing in Lemma~\ref{lm:push}, then we shall define the commutator
\[
 [D^m_s,\mathcal{L}]u:=D^m_s\mathcal{L}u-\mathcal{L}D^m_su.
\]
We formulate the following model problem, which is then used to derive the high-order energy identities.
Corresponding to the equations~(\ref{eq:tempu}) -~(\ref{eq:neumannu}) we 
analyze the equations:
\be\label{eq:prva}
\U_t-a_{ij}\U_{ij}-b_i\U_i=\mathcal A\quad\mbox{in}\,\,\,\,\Omega,
\ee 
\be
\U\circ\phi=-(n-1)\chi-\Delta_g\chi+\mathcal B\quad\mbox{on}\,\,\,\,\S,
\ee
\be
[\partial_n\U]^+_-\circ\phi=\frac{r^2}{|g|}\big(-\frac{\chi_tr}{|g|}+\alpha(t)\cdot\vec{n}_{\Gamma}\circ\phi
+\mathcal C\big)\quad\mbox{on}\,\,\,\,\S.
\ee
\begin{lemma}
\label{lm:energy}
The following energy identities hold:
\begin{enumerate}
\item[$(1)$]
\be\label{eq:en1}
\begin{array}{l}
\displaystyle
 \frac{1}{2}\partial_t\int_{\Omega}\rho\U^2
+\frac{1}{2}\partial_t\int_{\S}\big\{|\nabla_g\chi|^2-(n-1)\chi^2\big\}
+\int_{\Omega}\rho\nabla\U^tA\nabla\U\\
\displaystyle
\qquad=\int_{\Omega}P+\int_{\S}Q+\int_{\partial\Omega}\U\U_n;
\end{array}
\ee
\item[$(2)$]
\be\label{eq:en2}
\begin{array}{l}
\displaystyle
\int_{\Omega}\rho\U_t^2+\frac{1}{2}\partial_t\int_{\Omega}\rho|D\U|^2
+\partial_t\int_{\Omega}\rho\nabla\U^tA\nabla\U
+\int_{\Omega}\rho a_{ij}D_k\U_iD_k\U_j\\
\displaystyle
+\int_{\S}\big\{|\nabla_g\chi_t|^2-(n-1)|\chi_t|^2\big\}
+\frac{1}{2}\partial_t\int_{\S}\big\{|\nabla_g^2\chi|^2-(n-1)|\nabla_g\chi|^2\big\}\\
\displaystyle
\qquad=\int_{\Omega}S+\int_{\S}T+\int_{\g\Omega}\big[-(\U_t+\Delta\U)\U_n+\nabla\U\cdot\nabla\U_n\big],
\end{array}
\ee
where 
\be\label{eq:PQ}
P:=-(a_{ij}\rho)_{x^j}\U_i\U+(\mathcal{A}+b_i\U_i)\rho\U;
\quad
Q:=\underline{\frac{|g|}{r}\alpha(t)\cdot\vec{n}_{\Gamma}\circ\phi\,\U}
+\chi_t\mathcal{B}-\frac{|g|}{r}\mathcal{C}\U,
\ee
\be\label{eq:S}
\begin{array}{l}
\displaystyle
S:=\big(\mathcal A+b_i\U_i\big)\big(\U_t-\rho D_iD_i\U\big)
+\frac{1}{2}\rho_t|D\U|^2
+(\rho a_{ij})_t\U_i\U_j\\
\displaystyle
\qquad-(a_{ij}\rho)_j\U_t\U_i
+(a_{ij}\rho)_j\U_iD_kD_k\U-D_k(a_{ij}\rho)\U_iD_k\U_j\\
\displaystyle
\qquad-\underline{\rho D_i\U\U_t\nu(\mu)}
-\underline{a_{ij}\rho\U_iD_k\U_j\nu(\mu)};
\end{array}
\ee
\be\label{eq:T}
T:=-\chi_t\mathcal{B}_t+\frac{|g|}{r}\mathcal{C}\U_t
-\chi_t\Delta_g\mathcal{B}+\frac{|g|}{r}\mathcal{C}\Delta_g\U
+\underline{\frac{r}{|g|}\alpha(t)\cdot\vec{n}_{\Gamma}\circ\phi\,\U_t}
-\underline{\frac{r}{|g|}\alpha(t)\cdot\vec{n}_{\Gamma}\circ\phi\,\Delta_g\U}.
\ee
\end{enumerate}
\end{lemma}
\prf
We multiply~(\ref{eq:prva}) by $\rho\U$ and integrate over $\Omega$. Integrating by parts, we obtain
\[
\begin{array}{l}
\displaystyle
 \frac{1}{2}\partial_t\int_{\Omega}\rho\U^2+\int_{\Omega}\rho\nabla\U^tA\nabla\U
-\int_{\S}a_{ij}[\U_i]n^j\rho\U\\
\displaystyle
\qquad=\frac{1}{2}\int_{\Omega}\rho_t\U^2
-\int_{\Omega}(a_{ij}\rho)_{x^j}\U_i\U+\int_{\Omega}\big(\mathcal A+b_i\U_i\big)\rho\U
+\int_{\g\Omega}\U\U_n.
\end{array}
\]
Note that
\be
\begin{array}{l}
\displaystyle
 -\int_{\S}a_{ij}[\U_i]n^j\rho\U=-\int_{\S}a_{ij}n^in^jr[\U_n]\U \\
\displaystyle
=-\int_{\S}\frac{|g|^2}{r^4}r\frac{r^2}{|g|}\big(-\frac{\chi_tr}{|g|}+\alpha(t)\cdot\vec{n}_{\Gamma}\circ\phi+\mathcal C\big)
\big(-(n-1)\chi-\Delta_g\chi+\mathcal B\big)\\
\displaystyle
=\frac{1}{2}\partial_t\int_{\S}\big\{|\nabla_g\chi|^2-(n-1)\chi^2\big\}
-\underline{\int_{\S}\frac{|g|}{r}\alpha(t)\cdot\vec{n}_{\Gamma}\circ\phi\,\U}
+\int_{\S}-\chi_t\mathcal{B}+\frac{|g|}{r}\mathcal{C}\U.
\end{array}
\ee
This proves the first identity.
As to the second identity, we multiply~(\ref{eq:prva}) by $\rho\big(\U_t-D_iD_i\U\big)$, 
integrate over $\Omega$ and use integration by parts:
\be
\begin{array}{l}
\displaystyle
\int_{\Omega}\rho\U_t^2+\frac{1}{2}\partial_t\int_{\Omega}\rho|D\U|^2
+\partial_t\int_{\Omega}\rho\nabla\U^tA\nabla\U
+\int_{\Omega}\rho a_{ij}D_k\U_iD_k\U_j\\
\displaystyle
-\int_{\S}a_{ij}n^in^jr[\U_n]\U_t+\int_{\S}a_{ij}n^in^jr[\U_n]^+_-\Delta_g\U
=\int_{\Omega}\big(\mathcal A+b_i\U_i\big)\big(\U_t-\rho D_iD_i\U\big)\\
\displaystyle
\quad+\frac{1}{2}\int_{\Omega}\rho_t|D\U|^2-\underline{\int_{\Omega}\rho D_i\U\U_t\nu(\mu)}
+\int_{\Omega}(\rho a_{ij})_t\U_i\U_j-\int_{\Omega}(a_{ij}\rho)_j\U_t\U_i\\
\displaystyle
+\int_{\Omega}(a_{ij}\rho)_j\U_iD_kD_k\U-\int_{\Omega}D_k(a_{ij}\rho)\U_iD_k\U_j
-\underline{\int_{\Omega}a_{ij}\rho\U_iD_k\U_j\nu(\mu)}\\
\displaystyle
\quad+\int_{\g\Omega}\big[-(\U_t+\Delta\U)\U_n+\nabla\U\cdot\nabla\U_n\big].
\end{array}
\ee 
Let us extract the energy contribution from the integrals over $\S$ in the above identity.
\[
\begin{array}{l}
\displaystyle
-\int_{\S}a_{ij}n^in^jr[\U_n]\U_t
=-\int_{\S}\frac{|g|^2}{r^4}r\frac{r^2}{|g|}
\big(-\frac{\chi_tr}{|g|}+\alpha(t)\cdot\vec{n}_{\Gamma}\circ\phi+\mathcal C\big)
\big(-(n-1)\chi_t-\Delta_g\chi_t+\mathcal B_t\big)\\
\displaystyle
=\int_{\S}\big\{|\nabla_g\chi_t|^2-(n-1)|\chi_t|^2\big\}
-\underline{\int_{\S}\frac{r}{|g|}\alpha(t)\cdot\vec{n}_{\Gamma}\circ\phi\,\U_t}
+\int_{\S}\chi_t\mathcal{B}_t-\frac{|g|}{r}\mathcal{C}\U_t.
\end{array}
\]
Similarly,
\[
\begin{array}{l}
\displaystyle
\int_{\S}a_{ij}n^in^jr[\U_n]^+_-\Delta_{\mathcal B}\U
=\int_{\S}\frac{|g|^2}{r^4}r\frac{r^2}{|g|}
\big(-\frac{\chi_tr}{|g|}+\alpha(t)\cdot\vec{n}_{\Gamma}\circ\phi+\mathcal C\big)
\big(-(n-1)\Delta_g\chi-\Delta_g^2\chi+\Delta_g\mathcal B\big)\\
\displaystyle
=\frac{1}{2}\partial_t\int_{\S}\big\{|\nabla_g^2\chi|^2-(n-1)|\nabla_g\chi|^2\big\}
+\underline{\int_{\S}\frac{r}{|g|}\alpha(t)\cdot\vec{n}_{\Gamma}\circ\phi\,\Delta_g\U}
+\int_{\S}\chi_t\Delta_g\mathcal{B}-\frac{|g|}{r}\mathcal{C}\Delta_g\U.
\end{array}
\]
From the two previous identities we obtain the expression~(\ref{eq:T}) for the error term $T$.
\prfe
\subsection{The energy}\label{su:energy}
We introduce a parameter $\gamma>1$ to be fixed later. Motivated by Lemma~\ref{lm:energy},
for given $\U, \chi$ we define the auxiliary energy
\be\label{eq:energyaux}
\begin{array}{l}
\displaystyle
\E(\U,\chi):= \frac{\gamma}{2}\int_{\Omega}\rho\U^2
+\frac{1}{2}\int_{\Omega}\rho|D\U|^2+\int_{\Omega}\rho\nabla\U^tA\nabla\U\\
\displaystyle
\qquad+\frac{\gamma}{2}\int_{\S}\big\{|\nabla_g\chi|^2-(n-1)\chi^2\big\}
+\frac{1}{2}\int_{\S}\big\{|\nabla_g^2\chi|^2-(n-1)|\nabla_g\chi|^2\big\}
\end{array}
\ee
and the auxiliary dissipation
\be\label{eq:dissipationaux}
\D(\U,\chi):=\gamma\int_{\Omega}\rho\nabla\U^tA\nabla\U+\int_{\Omega}\rho\U_t^2
+\int_{\Omega}\rho a_{ij}D_k\U_iD_k\U_j
+\int_{\S}\big\{|\nabla_g\chi_t|^2-(n-1)|\chi_t|^2\big\}.
\ee
Summing the two previous expression and using Lemma~\ref{lm:energy}, we obtain
\be\label{eq:mainidentity}
\begin{array}{l}
\displaystyle
\E(\U,\chi)(t)+\int_0^t\D(\U,\chi)(s)\,ds=\E(\U,\chi)(0)+\int_0^t\int_{\Omega}\big\{\gamma P+S\big\}
+\int_0^t\int_{\S}\big\{\gamma Q+T\big\}\\
\displaystyle
\qquad+\int_0^t\int_{\partial\Omega}\big[\U_n(\gamma\U-\U_t-\Delta\U)+\nabla\U\cdot\nabla\U_n\big].
\end{array}
\ee
We define the total energy and dissipation by setting
\be\label{eq:en}
\E(w,f)=\sum_{|m|+2s\leq 2N}\E(D^m_sw,D^m_sf),
\ee
\be\label{eq:di}
\D(w,f)=\sum_{|m|+2s\leq 2N}\D(D^m_sw,D^m_sf).
\ee
We shall often write $\E(t):=\E(w(t),f(t))$ and $\D(t):=\D(w(t),f(t))$.
\section{A-priori estimates}\label{se:apriori}
\begin{lemma}[Poincar\'e type estimate]
\label{lm:poincare}
There exists a positive constant $\beta$ such that on the time interval of existence of 
a solution to the Stefan problem, the following estimate holds:
\[
||f||_{L^2}+||w||_{L^2(\Omega)}\leq \beta||\nabla w||_{L^2(\Omega)}.
\]
\end{lemma}
\prf
To prove the lemma, it is much more instructive to work with the solution $u$ on the moving domain
$\Omega(t)$ and later convert it into a result for $w$.  
Recall the expansion:
$
u\circ\phi=-(n-1)f-\Delta_gf+N(f).
$
Fix a spherical harmonic $s_{k,i}$ of degree $k$, where $k\geq2$ 
and $1\leq i\leq N(k)$ (recall the notation introduced in the introduction).
Multiply both sides of the
above relation by $s_{k,i}$ and integrate over $\S$. Observing that
$\int_{\S}(-(n-1)f-\Delta_gf)s_{k,i}=(n-1)(k^2-1)f_{k,i}$, we get
\be\label{eq:decaya}
f_{k,i}=\frac{1}{(n-1)(k^2-1)}\int_{\S}u\circ\phi\,s_{k,i}-\frac{1}{(n-1)(k^2-1)}\int_{\S}N(f)\,s_{k,i}.
\ee
Note that
\beas
\Big|\int_{\S}u\circ\phi\,s_{k,i}\Big|&=&
\Big|\int_{\S}\big(u\circ\phi-\frac{1}{|\S|}\int_{\S}u\circ\phi\big)\,s_{k,i}\Big|\\
&\leq& C||u\circ\phi-\frac{1}{|\S|}\int_{\S}u\circ\phi||_{L^2}\leq C||\nabla u||_{L^2(\Omega^{\pm})},
\eeas
where we used the Sobolev inequality in the last estimate above.
The last inequality, together with~(\ref{eq:decaya}), immediately implies 
\be\label{eq:decayak}
|f_{k,i}|\leq C\sqrt{\D}+C\Big|\int_{\S}N(f)\,s_{k,i}\,d\xi\Big|,\,\,\,\,k\geq2
\ee
In order to estimate $\int_{\S}f$, we observe that
\[
u\circ\phi r^{n-2}|g|=-(n-1)f-\Delta_gf+q(f),
\]
where $q$ stands for the nonlinear remainder with a leading order quadratic term. Integrating
the above equation over $\S$, we find
\be\label{eq:decay3}
\int_{\Gamma(t)}u=-(n-1)\int_{\S}f+\int_{\S}q(f),
\ee
Multiplying the conservation law
$
\int_{\Omega}u+\int_{\S}f+\sum_{k=2}^n{n\choose k}\int_{\S}f^k/n=0
$
by $\frac{1}{|\Omega|}$ 
and~(\ref{eq:decay3}) by $\frac{1}{|\S|}$ and subtracting the two equations, we obtain
\be
\label{eq:decay4}
\big(\frac{1}{|\Omega|}-\frac{n-1}{|\Gamma|}\big)\int_{\S}f=
-\big(\frac{1}{|\Omega|}\int_{\Omega}u-\frac{1}{|\Gamma(t)|}\int_{\Gamma(t)}u\big)
-\frac{1}{n|\Omega|}\sum_{k=2}^n{n\choose k}\int_{\S}f^k+\frac{1}{|\S|}\int_{\S}q(f).
\ee
Note that $\zeta-\big(\frac{1}{|\Omega|}-\frac{n-1}{|\Gamma(t)|}\big)=\frac{n-1}{|\Gamma(t)|}-\frac{n-1}{\S}$ and
hence $|\zeta-\big(\frac{1}{|\Omega|}-\frac{n-1}{|\Gamma(t)|}\big)|\leq C||f||_{H^1}$, 
which for $||f||_{H^1}$ small enough
implies $|\frac{1}{|\Omega|}-\frac{n-1}{|\Gamma(t)|}|\geq\zeta/2>0$.
Hence, upon dividing~(\ref{eq:decay4})
by $K_1:=\frac{1}{|\S|}\big(\frac{1}{|\Omega|}-\frac{n-1}{|\Gamma(t)|}\big)$, we conclude
\be\label{eq:decay5}
\int_{\S}f=-\frac{1}{K_1}\Big(\frac{1}{|\Omega|}\int_{\Omega}u-\frac{1}{|\Gamma(t)|}\int_{\Gamma(t)}u\Big)
-\frac{1}{nK_1|\Omega|}\sum_{k=2}^n{n\choose k}\int_{\S}f^k
+\frac{1}{K_1|\Gamma(t)|}\int_{\S}q(f).
\ee
From the mean value theorem, we deduce
\[
\Big|\frac{1}{|\Omega|}\int_{\Omega}u-\frac{1}{|\Gamma(t)|}\int_{\Gamma}u\Big|
\leq C||\nabla u||_{L^2(\Omega^{\pm})}.
\]
Thus, from~(\ref{eq:decay5}) and the previous inequality, we have
\be\label{eq:decay0}
|\int_{\S}f|\leq C||\nabla u||_{L^2(\Omega^{\pm})}+C||f||_{L^2}^2
\ee
Summing~(\ref{eq:decay0}) and~(\ref{eq:decayak}) and keeping in mind that $\q_1f=0$, we obtain
\[
\begin{array}{l}
\displaystyle
|\int_{\S}f|+\sum_{k=1}^{\infty}\sum_{i=1}^{N(k)}|f_{k,i}|
\leq C||\nabla u||_{L^2(\Omega)}+C||f||_{L^2}^2
+C\sum_{k=2}^{\infty}\sum_i\Big\{\Big|\int_{\S}N(f)s_{k,i}\Big|\\
\displaystyle
\leq C||\nabla u||_{L^2(\Omega)}+C||f||_{L^2}^2+C||N(f)||_{L^2}
\leq C||\nabla u||_{L^2(\Omega)}+C||f||_{L^2}^2.
\end{array}
\] 
Smallness of $||f||_{L^2}$ then implies
\[
\sum_{k=0}^{\infty}\sum_{i=1}^{N(k)}|f_{k,i}|
\leq C||\nabla u||_{L^2(\Omega(t)}.
\]
Changing variables ($x\to\Theta(x)$, i.e. $u\to w$) and using the smallness of $||\nabla_gf||_{L^2}$
we get the bound $||\nabla u||_{L^2(\Omega^{\pm})}\leq C||\nabla w||_{L^2(\Omega)}$. Combined with 
previous estimate we conclude
\be
\label{eq:decay8}
||f||_{L^2}\leq C||\nabla w||_{L^2(\Omega(t))}
\ee
Finally, from the previous inequality and the conservation law~(\ref{eq:mass}), we immediately
deduce 
\[
\Big|\int_{\Omega}u\Big|\leq C||\nabla u||_{L^2(\Omega(t))}\leq C||\nabla w||_{L^2(\Omega(t))}.
\]
Changing variables and using the smallness of $||\nabla_gf||_{L^2}$, we obtain
\be
\label{eq:decay9}
\Big|\int_{\Omega}w\Big|\leq C||\nabla w||_{L^2(\Omega^{\pm})}
\ee
and this finishes the proof of the lemma.
\prfe
In the new coordinates, the mass conservation law~(\ref{eq:mass}) takes a different form. The 
following lemma expresses this conservation in a way that will be useful for proving the positive
definitiveness of the energy expressions $\E$ and $\D$.
\begin{lemma}[Mass conservation law]
\label{lm:consfixed}
The following identity holds:
\be\label{eq:consfixed}
\int_{\Omega}\rho w+\int_{\S}f=-\int_{\S}\sum_{k=2}^n{n \choose k}\frac{f^k}{n}+\int_{\Omega}wg(\nabla\rho, f),
\ee
where $g$ is a bounded smooth function with $g({\bf 0},0)=0$.
\end{lemma}
\prf
Integrating the equation~(\ref{eq:tempu}) on
a moving domain $\Omega\setminus\Gamma(t)$, using the Stokes formula and the boundary 
condition~(\ref{eq:neumannu}), we obtain:
\[
\g_t\int_{\Omega(t)}u+\g_t\int_{\S}\frac{(1+f)^n}{n}=0.
\]
The assumed initial condition~(\ref{eq:initialmassu}) implies
\[
\int_{\Omega(t)}u+\sum_{k=1}^n{n\choose k}f^k=0.
\] 
To express the integral over $\Omega(t)$ above as an integral over $\Omega\setminus\S$, 
it remains to understand the Jacobian $|\mbox{det}D\Theta^{-1}|$, 
where $\Theta:\Omega(t)\to\Omega\setminus\S$ is the change of variables map
defined in the line before~(\ref{eq:change1}) and $\Theta^{-1}:\bar{x}\to \a(t)+\rho(t,\bar{x})\bar{x}$.
Thus, $(D\Theta^{-1})_{ij}=\rho\delta_{ij}+\rho_{\bar{x}^i}\bar{x}^j$, where $\delta_{ij}$ denotes
the Kronecker delta. This implies $\mbox{det}D\Theta^{-1}=\rho^n+q_1(\nabla\rho)
=\rho+(\rho-1)q_2(\rho)+q_1(\nabla\rho)$, for some polynomials $q_1$ and $q_2$ 
such that $q_i(0)=0$, $i=1,2$. From here, we easily infer the lemma. Since $\rho-1=f$ close to
$\S$ and $\rho$ is smooth and equal to $0$ close to the boundary $\g\Omega$, the claim of the lemma
follows.
\prfe
\begin{lemma}[Positivity of the energy]\label{lm:positive}
Under the smallness assumption on $||D^m_sf||_{L^2}+||D^m_sw||_{L^2}$, for $|m|+2s\leq N$, 
the energy quantities $\E$ and $\D$ are positive definite.
\end{lemma}
\prf
For a given function $\omega:\S\to\R$, let us abbreviate 
$Z(\omega)=\int_{\S}|\nabla_g\omega|^2-(n-1)\omega^2$.
If $\chi=D^m_sf$, with $|\mu|\geq1$, then, by Wirtinger's inequality, we immediately see
\[
Z(\chi)\geq C||\q_{2+}\chi||_{L^2}^2
\]
and analogous inequalities hold for $Z(\chi_t)$ and $Z(\nabla\chi)$. If however, $\chi$ is 
of the form $\chi=\partial_{t^s}f$, then we must exploit the conservation law~(\ref{eq:consfixed}).
Note that
\be
\label{eq:pos1}
\int_{\Omega}\rho w_{t^s}^2+Z(f_{t^s})
=\frac{1}{|\Omega|_{\rho}}\big(\int_{\Omega}\rho w_{t^s}\big)^2
+\int_{\Omega}\rho(w_{t^s}-\frac{1}{|\Omega|_{\rho}}\int_{\Omega}\rho w_{t^s})^2
-\frac{(n-1)}{|\S|}\big(\int_{\S}f_{t^s}\big)^2+Z(\q f_{t^s}),
\ee
where $|\Omega|_{\rho}=\int_{\Omega}\rho\,dx$.
From Lemma~\ref{lm:consfixed}, we obtain
\[
\int_{\Omega}\rho w_{t^s}+\int_{\S}f_{t^s}=G(w,f),
\]
where 
$
G(w,f)=\int_{\Omega}\big\{-\partial_{t^s}\big(\rho w\big)+\rho w_{t^s}\big\}
-\int_{\S}\sum_{k=2}^n{n \choose k}f^k_{t^s}/n+\int_{\Omega}\partial_{t^s}(wg(\nabla\rho,f))
$.
Thus, 
\[
\begin{array}{l}
\displaystyle
\frac{1}{|\Omega|_{\rho}}\big(\int_{\Omega}\rho w_{t^s}\big)^2
-\frac{(n-1)}{|\S|}\big(\int_{\S}f_{t^s}\big)^2\\
\displaystyle
\quad=\big(\frac{1}{|\Omega|_{\rho}}-\frac{n-1}{|\S|}\big)\big(\int_{\Omega}\rho w_{t^s}\big)^2
-G(w,f)\big(G(w,f)-2\int_{\Omega}\rho w_{t^s}\big).
\end{array}
\]
Note however
\[
\frac{1}{|\Omega|_{\rho}}-\frac{n-1}{|\S|}=\zeta+\frac{1}{|\Omega|}-\frac{1}{|\Omega|_{\rho}}
=\zeta+h(f),
\]
where $h$ is a smooth bounded function, $h(0)=0$. Going back to~(\ref{eq:pos1}), we conclude
\[
\int_{\Omega}\rho w_{t^s}^2+Z(f_{t^s})
=\zeta\big(\int_{\Omega}\rho w_{t^s}\big)^2
+\int_{\Omega}\rho(w_{t^s}-\frac{1}{|\Omega|_{\rho}}\int_{\Omega}\rho w_{t^s})^2
+Z(\q f_{t^s})+R(w_{t^s},f_{t^s})
\]
where the remainder $R(w_{t^s},f_{t^s})$ is a cubic nonlinearity given by
\[
R(w_{t^s},f_{t^s}):=-G(w,f)\big(G(w,f)-2\int_{\Omega}\rho w_{t^s}\big)+h(f)\big(\int_{\Omega}\rho w_{t^s}\big)^2.
\]
Under a smallness assumption on $w$ and $f$ it is easy to see that
$
\big|R(w_{t^s},f_{t^s})\big|\leq \lambda\big(||w_{t^s}||_{L^2(\Omega)}^2+||f_{t^s}||_{L^2}^2\big),
$
for a small constant $\lambda$ and thus there exists a positive constant $M\geq1$ such that
\[
\frac{1}{M}\big(||w_{t^s}||_{L^2(\Omega)}^2+||f_{t^s}||_{L^2}^2\big)
\leq\int_{\Omega}\rho w_{t^s}^2+Z(f_{t^s})\leq M\big(||w_{t^s}||_{L^2(\Omega)}^2+||f_{t^s}||_{L^2}^2\big).
\]
This concludes the proof of the lemma.
\prfe
The following lemma states that the norms defined by $\E$ and $\D$ are equivalent to the
corresponding parabolic Sobolev norms. The main new ingredient is that {\em all} the derivatives of $w$
are bounded. In particular, derivatives of $w$ in the normal direction are bounded by the energy quantities $\E$ and
$\D$, although they are not a-priori contained in their definitions~(\ref{eq:en}) and~(\ref{eq:di}).
\begin{lemma}\label{lm:normequiv}
The following norms are equivalent:
\[
\E(w,f)(t)\approx
\sum_{k=0}^N\big\{||w||_{W^{N-k,\infty}W^{2k,2}([0,t]\times\Omega)}^2
+||\nabla_gf||_{W^{N-k,\infty}W^{2k+1,2}([0,t]\times\K^{n-1})}^2\big\}.
\]
\[
\begin{array}{l}
\displaystyle
\D(w,f)(t)\approx\sum_{k=0}^N\big\{||\nabla w||_{W^{N-k,2}W^{2k+1,2}([0,t]\times\Omega)}^2
+||w_t||_{W^{N-k,2}W^{2k,2}([0,t]\times\Omega)}^2\\
\displaystyle
\qquad\qquad\qquad+||f_t||_{W^{N-k,2}W^{2k+1,2}([0,t]\times\Omega)}^2\big\},
\end{array}
\]
where $A\approx B$ means that there is a constant $c>1$ such that $\frac{A}{c}\leq B\leq cA$.
\end{lemma}
\prf
We only sketch the proof as it is analogous to the proof of Lemma 3.5 in~\cite{HaGu1}. 
What we need to prove is that for any triple of indices $(\mu,r,s)$ satisfying
$|\mu|+r+2s\leq2N$ we have
\[
\|D_{\xi}^{\mu}\g_n^r\g_{t^s}w\|_{L^2(\Omega)}^2\leq C\E(w,f).
\]
The claim is evidently true for $r=0,1$ following from the definition of $\E$. We then proceed by induction 
on the number of normal derivatives $r$. Inductive step uses the basic equation~(\ref{eq:tempw}). It
allows to express $w_{nn}$ as the sum of the multiples of $w_t$, $D^{ij}w$ and $D^iw_n$ in the
region of $\Omega$ close to the unit sphere $\S$, thus allowing to complete the proof in the case $r=2$. 
By successively differentiating the equation~(\ref{eq:tempw}) and expressing the term with the highest number
of normal derivatives in terms of terms with either {\it less} normal derivatives or {\it purely tangential}
derivatives, we inductively complete the argument.
\prfe
\section{Local existence}\label{se:local}
\begin{theorem}\label{th:local}
There exist small positive constants $E_1$, $m_1$ such that for any
$E_0\leq E_1$, $m_0\leq m_1$ and positive
constants $T^*$, $C^*$, such that if
\[
\E(w_0,f_0)\leq E_0;\quad|\a(0)|\leq m_0,
\]
then there exists a unique solution to the Stefan 
problem~(\ref{eq:tempu}) -~(\ref{eq:initialmassu}) on the time
interval $[0,T^*]$ such that for any $0\leq s\leq t<T^*$:
\be\label{eq:apriori}
\sup_{s\leq\tau\leq t}\E(\tau)+\int_s^t\D(\tau)\,d\tau
\leq\E(s)+(\frac{1}{4}+\tilde{C}\sup_{s\leq\tau\leq t}\sqrt{\E(\tau)})\int_0^t\D(\tau)\,d\tau;
\ee
\be\label{eq:apriori2}
\sup_{0\leq\tau\leq T^*}|a(\tau)|\leq2m_0.
\ee
Moreover, $\E$ is continuous on $[0,T^*[$ and
$\sup_{s\leq\tau\leq t}\E(\tau)+\frac{1}{2}\int_s^t\D(\tau)\,d\tau\leq E_0$.
\end{theorem}
\prf We first sketch the proof of the a-priori bound~(\ref{eq:apriori}) assuming that the solution
already exists. We explain in detail how to 
bound the hardest quadratic-in-order terms. For all the other (cubic) terms
we use rather standard energy estimates, that can be found
in the corresponding local existence proofs in~\cite{HaGu1} (or~\cite{HaGu2}).
To prove the estimate~(\ref{eq:apriori}), we must bound the right-hand side
of the energy identity~(\ref{eq:mainidentity}) where we plug in
$\U=D^m_sw$ and $\chi=D^m_sf$ for all $|m|+2s\leq2N$. 
All the terms in the definitions
of $P$, $Q$, $S$, and $T$ (cf.~(\ref{eq:PQ}) -~(\ref{eq:T})) are trilinear
except for the underlined terms, which are 
only quadratic in the order of nonlinearity and it is a-priori not clear how 
to bound them. 
\par
There are two types of the quadratic error terms: the first type arises
for technical reasons due to the introduction of the cut-off function $\mu$
while fixing the domain in Section~\ref{se:fix}.
These are the two underlined terms in the expression~(\ref{eq:S}) for $S$.
For $\U=D^m_sw$ and $\chi=D^m_sf$, $|m|+2s\leq2N$,
by the Young's inequality, we have
\be\label{eq:gamma1}
\begin{array}{l}
\displaystyle
\Big|\int_0^t\int_{\Omega}\rho D_i\U\U_t\nu(\mu)\Big|
\leq\frac{1}{4}\int_0^t\int_{\Omega}\rho\U_t^2
+||\nu(\mu)||_{L^{\infty}(\Omega)}^2\int_0^t\int_{\Omega}\rho|D\U|^2\\
\displaystyle
\qquad\leq\frac{1}{4}\int_0^t\int_{\Omega}\rho\U_t^2+C^*\int_0^t\int_{\Omega}\rho\nabla\U^tA\nabla\U.
\end{array}
\ee 
Similarly we estimate the integral of the second underlined term in 
the expression~(\ref{eq:S}):
\be\label{eq:gamma2}
\Big|\int_0^t\int_{\Omega}a_{ij}\rho\U_iD_k\U_j\nu(\mu)\Big|
\leq\frac{1}{4}\int_0^t\int_{\Omega}a_{ij}\rho D_k\U_jD_k\U_i
+C^*\int_0^t\int_{\Omega}\rho\nabla\U^tA\nabla\U,
\ee
where we possibly enlarge the constant $C^*$. The size of the constant $C^*$ dictates the choice 
of $\gamma$ in the definition of the energy (cf. Section~\ref{su:energy}): 
set $\gamma=3C^*$. We can thus absorb the right-most terms in the
estimates~(\ref{eq:gamma1}) and~(\ref{eq:gamma2}) into the $\gamma$-dependent term in the 
definition~(\ref{eq:dissipationaux}) of $\D$.
\par
The second type of quadratic error terms is completely {\em intrinsic} and arises due to the
{\em presence of the moving center coordinate $\a(t)$} in the parametrization of the moving surface
$\Gamma(t)$. 
They are the three terms, underlined in 
the expressions~(\ref{eq:PQ}) and~(\ref{eq:T}) for $Q$ and $T$ respectively. 
Note that, if $|m|\geq1$, then $\alpha(t)=D^m_s\a(t)=0$, since $\a$ depends only on $t$.
Thus, we are only concerned with estimating expressions containing $\alpha(t)$ of the
form $\alpha(t)=\partial_t^{k+1}\a(t)=\a^{(k+1)}(t)$, where $0\leq k\leq N-1$. The underlined expression
in the formula for $Q$ thus takes the form
\[
\int_0^t\int_{\S}\frac{r}{|g|}\a^{(k+1)}(t)\cdot\vec{n}_{\Gamma}\circ\phi\,\partial_t^kw.
\]
The difficulty is immediately clear: since the expression $\frac{r}{|g|}\vec{n}_{\Gamma}\circ\phi$ is a quantity
of order $1$ for small $f$, the whole integral is only of a quadratic order and it is hence 
unclear how to bound it by $\sqrt{\E}\int_0^t\D$. In the local coordinates on the sphere: 
\[
\vec{n}_{\Gamma}\circ\phi=\frac{r}{|g|}\xi-\frac{\nabla_gf}{|g|},
\]
where we keep in mind that $\xi$ is the unit normal of the unit sphere $\S$. From here,
\[
\frac{r}{|g|}\vec{n}_{\Gamma}\circ\phi
=\xi+\frac{r^2-|g|^2}{|g|^2}\xi-\frac{r\nabla_gf}{|g|^2}.
\] 
Using this identity,
\[
\int_0^t\int_{\S}\frac{r}{|g|}\a^{(k+1)}(t)\cdot\vec{n}_{\Gamma}\circ\phi\,\partial_t^kw
=\int_0^t\int_{\S}\a^{(k+1)}(t)\cdot\big[\xi+\frac{r^2-|g|^2}{|g|^2}\xi-\frac{r\nabla_gf}{|g|^2}\big]
\,\partial_t^k(-(n-1)f-\Delta_gf+N(f))
\]
The key observation is $\xi\in\mbox{Null}(\Delta_g+(n-1)\mathcal I)$ where
$\mathcal{I}$ stands for the identity operator. Hence
\[
\int_{\S}\a^{(k+1)}(t)\cdot\xi\partial_{t^k}(\Delta_g+(n-1)\mathcal I)f=
\a^{(k+1)}(t)\cdot\int_{\S}(\Delta_g+(n-1)\mathcal I)\xi\partial_{t^k}f=0.
\]
Using the previous two identites, we eliminate the purely quadratic contribution:
\[
\begin{array}{l}
\displaystyle
\int_0^t\int_{\S}\a^{(k+1)}(t)\cdot\big[\xi+\frac{r^2-|g|^2}{|g|^2}\xi-\frac{r\nabla_gf}{|g|^2}\big]
\,\partial_{t^k}(-(n-1)f-\Delta_gf+N(f))\\
\displaystyle
=\int_0^t\int_{\S}\a^{(k+1)}(t)\cdot\big[\frac{r^2-|g|^2}{|g|^2}\xi-\frac{r\nabla_gf}{|g|^2}\big]
\,(-(n-1)f_{t^k}-\Delta_gf_{t^k})\\
\displaystyle
\quad
+\int_0^t\int_{\S}\a^{(k+1)}(t)\cdot\big[\frac{r^2}{|g|^2}\xi-\frac{r\nabla_gf}{|g|^2}\big]
\,\partial_{t^k}(N(f))
\end{array}
\]
Recall that $N(f)$ is a {\it quadratic} nonlinearity and therefore the two integrands on the right-hand 
side above have a manifestly trilinear 
structure and are thus easy to estimate. Namely, applying the differential
operator $\partial_{t^k}$ to~(\ref{eq:orthAB}), it is easy to deduce the bound 
$\big|\a^{(k+1)}(t)\big|\leq C\sqrt{\D}$ if we knew that $|\a(\tau)|$ is bounded. 
To estimate $|\a(\tau)|$ we use~(\ref{eq:orthAB}). Recalling ${\bf p}={\bf x}-\a$, 
and assuming smallness of $\E$, we obtain
\[
|\a(T)|\leq|\a(0)|+T\sup_{0\leq\tau\leq T}|\dot{\a}(\tau)|
\leq|\a(0)|+CT\sqrt{\E}+CT\sup_{0\leq\tau\leq T}|\a(\tau)|\|u_t\|_{L^2(\Omega)}.
\]
Therefore
\[
\sup_{0\leq\tau\leq T}|\a(\tau)|(1-CT\sqrt{\E})\leq|\a(0)|+CT\sqrt{\E},
\]
thus implying $\sup_{0\leq\tau\leq T}|\a(\tau)|\leq2|\a(0)|\leq2m_0$ for appropriately small $\E$ and 
some finite $T$.
Integration by parts, Lemma~\ref{lm:poincare} and
standard energy estimates then imply
\be\label{eq:orth1}
\Big|\int_0^t\int_{\S}\frac{r}{|g|}\a^{(k+1)}(t)\cdot\vec{n}_{\Gamma}\circ\phi\,\partial_t^kw\Big|
\leq C\sqrt{\E}\int_0^t\D.
\ee
The same idea as in the proof of~(\ref{eq:orth1}) works for the remaining two underlined integrals
appearing in $\int_0^t\int_{\S}T$ (cf.~(\ref{eq:T})), so we finally conclude
(for $\alpha=\a^{(k+1)}(t)$ and $\U=w_{t^{k+1}}$)
\be
\Big|\int_0^t\int_{\S}\frac{r}{|g|}\alpha(t)\cdot\vec{n}_{\Gamma}\circ\phi\,\U_t\Big|
+\Big|\int_0^t\int_{\S}\frac{r}{|g|}\alpha(t)\cdot\vec{n}_{\Gamma}\circ\phi\,\Delta_g\U\Big|
\leq C\sqrt{\E}\int_0^t\D.
\ee
The construction of the solution $(w,f)$ follows identically the construction scheme
for the local solution from~\cite{HaGu1}. 
We briefly summarize the main steps:
we set-up an iteration scheme, which
generates a sequence of iterates $\{(w^m,f^m)\}_{m\in\N}$ solving a sequence
of linear  parabolic problems. 
As in~\cite{HaGu1}, such an
iteration is well defined, but it breaks the natural energy setting due to the
lack of exact cancellations in the presence of cross-terms. 
We design the elliptic regularization
\[
-\frac{f_tr}{|g|}-\epsilon\frac{\Delta_gf_t}{r^{n-2}|g|}=[w_n]^+_-\circ\phi.
\]
to overcome this difficulty. 
For a fixed $\epsilon$,
we use it to prove that $\{(w^m,f^m)\}_{m\in\N}$ is a Cauchy sequence in the energy space.
Upon passing to the limit $m\to\infty$ we obtain a solution existing on an a-priori
$\epsilon$-dependent time interval $[0,T^{\epsilon}]$. But, in the limit $m\to\infty$
the dangerous cross-terms vanish. Since the elliptic regularization honors the
energy structure, we obtain an energy bound analogous to~(\ref{eq:apriori})
with $\epsilon$-{\em independent} coefficients. By continuity, the solution exists on
an $\epsilon$-indpendent time interval $[0,T]$. Finally, we pass to the limit as 
$\epsilon\to0$ to obtain the solution of the original problem.
\prfe
\section{Global existence and asymptotic stability}\label{se:global}
The proof of the main result Theorem~\ref{th:main}, 
is an immediate consequence of the following theorem and 
Lemma~\ref{lm:normequiv} about the equivalence of norms.
\begin{theorem}\label{th:global}
There exist small positive constants $E$ and $m$, such that
if 
\[
\E(w_0,f_0)\leq E;\quad|\a(0)|\leq m,
\]
then there exists a unique global-in-time solution to the Stefan 
problem~(\ref{eq:tempu}) -~(\ref{eq:initialmassu}) converging
asymptotically to some steady state solution $\Sigma(\bar{\a})$, where
\[
\a(t)\to\bar{\a}\quad\mbox{as}\,\,\,\,t\to\infty.
\]
Moreover, there are constants $c_1, c_2>0$ such that the following exponential decay estimate holds:
\[
\E(w(t),f(t))\leq c_1e^{-c_2t}.
\]
\end{theorem}
Before we prove the theorem, we will prove an important auxiliary estimate,
necessary for the proof of Theorem~\ref{th:global}.
\begin{lemma}\label{lm:decay}
Let $[0,\mathfrak{t}]$ be an interval of existence of solution to the Stefan 
problem~(\ref{eq:tempu}) -~(\ref{eq:initialmassu}) for which the
estimates~(\ref{eq:apriori})  and~(\ref{eq:apriori2}) hold. 
Then there exist constants $\alpha,\delta\in\R_+$
such that
\[
\E(t)\leq\frac{4\beta\E(0)}{t}e^{-\frac{\alpha t}{2}};\quad t\in[0,\mathfrak{t}[
\]
and
\be\label{eq:decayAB}
|\dot{\a}(t)|\leq\frac{c_{\beta\delta}\sqrt{\E(0)}}{\sqrt{t}}e^{-\frac{\alpha t}{4}};
\quad t\in[0,\mathfrak{t}[,
\ee
where $\beta>0$ is given by Lemma~\ref{lm:poincare} and $c_{\beta\delta}:=2\sqrt{\beta}\delta$.
\end{lemma}
\prf
As in~\cite{Ma}, p. $135$,
we define $V(s):=\int_s^{\mathfrak{t}}\E(\tau)\,d\tau$. From~(\ref{eq:apriori}) and
Lemma~\ref{lm:poincare}, we conclude that there exists a constant $\alpha>0$ such
that
\be\label{eq:decay1}
\sup_{s\leq\tau\leq t}\E(\tau)+\alpha\int_s^t\E(\tau)\,d\tau\leq\E(s).
\ee
Thus $\alpha V(s)\leq\E(s)$. On the other hand,
$
V'(s)=-\E(s)\leq-\alpha V(s),
$
which, in turn, implies:
$
V(s)\leq V(0)e^{-\alpha s}.
$
Integrate~(\ref{eq:decay1}) over the interval $[t/2,t]$ with respect to $s$, to obtain:
\[
\frac{t}{2}\sup_{s\leq\tau\leq t}\E(\tau)\leq V(\frac{t}{2})\leq V(0)e^{-\frac{\alpha t}{2}},
\]
implying
\[
\E(t)\leq\frac{2V(0)}{t}e^{-\frac{\alpha t}{2}}.
\]
Note that for $E_0$ and $m_0$ small enough 
($\E(0)<\frac{1}{32\tilde{C}^2}$) in Theorem~\ref{th:local}, by continuity and the 
inequality~(\ref{eq:apriori}), we obtain 
the bound $\tilde{C}\sup_{0\leq s<\mathfrak{t}}\sqrt{\E(s)}\leq\frac{1}{4}$.
Plugging back into~(\ref{eq:apriori}) and absorbing the right-most term into LHS,
we obtain $\int_0^t\D(\tau)\,d\tau\leq 2\E(0)$ for all $t\in[0,\mathfrak{t}[$.  
Hence, by Lemma~\ref{lm:poincare} $V(0)=\int_0^{\mathfrak{t}}\E(\tau)\,d\tau
\leq\beta\int_0^{\mathfrak{t}}\D(\tau)\,d\tau\leq2\beta\E(0)$ and we obtain
\be\label{eq:decay}
\E(t)\leq\frac{4\beta\E(0)}{t}e^{-\frac{\alpha t}{2}},
\ee
thus proving the first claim of the lemma.
From~(\ref{eq:orthAB}) and Jensen's inequality, we conclude that for some constant
$\delta>0$ on the time interval $[0,\mathfrak{t}]$ the following
inequality holds:
\beas
|\dot{\a}|&\leq&\frac{\delta}{4}\sum_{k=1}^n||f^k||_{L^2}||f_t||_{L^2}
+\frac{\delta}{4}||w_t||_{L^2(\Omega)}
+\frac{\delta}{4}||\nabla w||_{L^2(\Omega)}\\
&\leq&\frac{\delta}{4}(2+\sum_{k=1}^{n-1}(2\E(0)^k))\sqrt{\E(w,f)}\\
&\leq&\delta\sqrt{\E(w,f)}.
\eeas
By~(\ref{eq:decay}) we obtain
\[
|\dot{\a}(t)|\leq\frac{2\sqrt{\beta}\delta\sqrt{\E(0)}}{\sqrt{t}}e^{-\frac{\alpha t}{4}}
=\frac{c_{\beta\delta}\sqrt{\E(0)}}{\sqrt{t}}e^{-\frac{\alpha t}{4}}
\]
for any $t\in[0,\mathfrak{t}]$ and this finishes the proof of the lemma.
\prfe
{\bf Proof of Theorem~\ref{th:global}.}
Let $m=\frac{m_1}{4}$ ($m_1$ is given by Theorem~\ref{th:local}) and $E$ be such that 
$c_{\beta\delta}c_{\alpha}\sqrt{E}\leq\frac{m_1}{4}$,
whereby $c_{\alpha}:=\max_{y\geq T^*}\frac{\sqrt{y}e^{-\frac{\alpha y}{4}}}{(1-e^{-\frac{\alpha T^*}{4}})}$.
Define 
\be
\mathcal{T}:=\sup_{\mathfrak{t}}
\Big\{
\sup_{0\leq\tau\leq\mathfrak{t}}|\a(\tau)|\leq2m\,\,\,
\&\,\,\, (\ref{eq:apriori})\,\,\,\mbox{holds for any}\,\,\,\,s,t\in[0,\mathfrak{t}]
\Big\}.
\ee
Assume $\mathcal{T}<\infty$. With $T^*>0$ given by Theorem~\ref{th:local},
choose a unique $L\in\N$ such that
\[
\mathcal{T}\in[LT^*, (L+1)T^*[.
\]
Integrating over $[kT^*, (k+1)T^*]$ with $1\leq k\leq L-1$, we
obtain
\beas
|\a((k+1)T^*)||&\leq&|\a(kT^*)|
+\int_{kT^*}^{(k+1)T^*}|\dot{\a}(\tau)|\,d\tau\\
&\leq&|\a(kT^*)|+c_{\beta\delta}\sqrt{\E(0)}\int_{kT^*}^{(k+1)T^*}
\frac{e^{-\frac{\alpha\tau}{4}}}{\sqrt{\tau}}\,d\tau\\
&\leq&|\a(kT^*)|+\frac{c_{\beta\delta}\sqrt{E(0)}}{\sqrt{kT^*}}\int_{kT^*}^{(k+1)T^*}
e^{-\frac{\alpha\tau}{4}}\,d\tau\\
&\leq&|\a(kT^*)|+\frac{c_{\beta\delta}\sqrt{\E(0)}}{\sqrt{kT^*}}T^*e^{-\frac{\alpha kT^*}{4}}\\
&\leq&|\a(kT^*)|+c_{\beta\delta}\sqrt{\E(0)}\sqrt{T^*}e^{-\frac{\alpha kT^*}{4}}.
\eeas
Summing over $k=1,\dots, L-1$, we obtain
\beas
|\a(LT^*)|&\leq&|\a(T^*)|+c_{\beta\delta}\sqrt{\E(0)}\sqrt{T^*}
\sum_{k=1}^{L-1} e^{-\frac{\alpha kT^*}{4}}
=|\a(T^*)|+\frac{c_{\beta\delta}\sqrt{\E(0)}\sqrt{T^*}}{1-e^{-\frac{\alpha T^*}{4}}}
\big(e^{-\frac{\alpha T^*}{4}}-e^{-\frac{\alpha LT^*}{4}}\big)\\
&\leq&2m+\frac{c_{\beta\delta}\sqrt{T^*}}{1-e^{-\frac{\alpha T^*}{4}}}e^{-\frac{\alpha T^*}{4}}\sqrt{E}
\leq 2m+c_{\beta\delta}c_{\alpha}\sqrt{E}\leq\frac{3}{4}m_1,
\eeas
where the last inequality follows from the choice of $m$ and $E$ above.
If we denote $T^1:=LT^*+\frac{\mathcal{T}-LT^*}{2}$ we obtain
\beas
|\a(T^1)|&\leq&|\a(LT^*)|
+\int_{LT^*}^{T^1}|\dot{\a}(\tau)|\,d\tau
\leq\frac{3m_1}{4}+c_{\beta\delta}\sqrt{\E(0)}\sqrt{T^*}e^{-\frac{\alpha LT^*}{4}}\\
&\leq&\frac{3m_1}{4}+c_{\beta\delta}c_{\alpha}\sqrt{E}\leq\frac{3m_1}{4}+\frac{m_1}{4}
=m_1.
\eeas
Since $\E(T^1)\leq E\leq E_1$ (with $E$ possibly smaller, cf.~(\ref{eq:decay})) 
and $|\a(T^1)|\leq m_1$, we can extend the solution $(w,f)$ starting at $t_0=T^1$
to the time interval $[0,T^1+T^*[$, by Theorem~\ref{th:local}
(local existence theorem).
Same theorem guarantees the estimate
\[
||\a||_{L^{\infty}[0,T^1+T^*]}\leq 2m_1
\]
as well as the validity of 
the energy estimate~(\ref{eq:apriori}) on the time-interval $[0,T^1+T^*]$. Since we
assumed $\mathcal{T}<\infty$, we obtain $T^1+T^*>\mathcal{T}$ and this, together with
the continuity of $\E(w,f)$ and $\a$ contradicts 
the maximality of $\mathcal{T}$. It is now evident that the decay 
estimates~(\ref{eq:decayAB}) and~(\ref{eq:decay}) from 
Lemma~\ref{lm:decay} hold for all $t\in[0,\infty[$. This implies the
decay claim of Theorem~\ref{th:global} as well as the existence
of an asymptotic center $\bar{\a}\in\Omega$ such that $\a(t)\to\bar{\a}$, since 
$\dot{\a}(t)\to0$ as $t\to\infty$.  
\prfe

\appendix

\section{Proof of Lemma~\ref{lm:push}}
Note that~(\ref{eq:change1}) implies 
\be\label{eq:simple}
\bar{x}_{x^i}=\pi e_i+\pi_{x^i}(x-\a).
\ee
From $\pi(x)=1/\rho(\bar{x})$, we obtain $\pi_{x^i}=-(\nabla\rho\cdot\bar{x}_{x^i})/\rho^2$,
which in turn, combined with~(\ref{eq:simple}), after an elementary calculation implies
\be\label{eq:simple2}
\pi_{x^i}=-\frac{\rho_{\bar{x}^i}}{\rho^2\Psi},
\ee
where
\[
 \Psi=\rho+\nabla\rho\cdot\bar{x}.
\]
After further differentiating~(\ref{eq:simple2}) with respect to $x^i$ and
using the relation~(\ref{eq:simple}), we arrive at
\[
\Delta\pi=\frac{-\Delta\rho}{\rho^3\Psi}
+\frac{\rho_{\bar{x}^i\bar{x}^j}\rho_{\bar{x}^i}\bar{x}^j(1-\rho^2)}{\rho^5\Psi^2}
+\frac{-\rho_{\bar{x}^i\bar{x}^j}\bar{x}^i\bar{x}^j|\nabla\rho|^2+2|\nabla\rho|^2(\Psi+\rho)}{\rho^3\Psi^3}.
\]
Similarly it is not hard to see that
$
\pi_t=(-\rho_t+\rho\nabla\rho\cdot\dot{\a})/\rho\Psi.
$
In order to evaluate the Laplacian in new coordinates, by~(\ref{eq:uchange}) we first write
\be\label{eq:uchange2}
u(t,x)=w(t,\pi(t,x)x).
\ee
Applying $\Delta_{x}$ to the left-hand side, we obtain
\be
\begin{array}{l}
\Delta u(t,x)=
\pi^2\Delta_{\bar{x}}w+\big(2\pi x^j\pi_{x^i}+x^ix^j|\nabla\pi|^2\big)
w_{\bar{x}^i\bar{x}^j}
+\big(2\pi_{x^i}+\Delta\pi   x^i\big)w_{\bar{x^i}}
\end{array}
\label{eq:laplacenew}
\ee
Applying $\partial_t$ to the left hand side of~(\ref{eq:uchange2}), we obtain
\be\label{eq:timenew}
u_t(t,x)=w_t+\pi_t\nabla_{\bar{x}}w\cdot x.
\ee
Since $(\partial_t-\Delta)u=0$, we conclude from~(\ref{eq:laplacenew}) and~(\ref{eq:timenew})
\[
(\partial_t-\Delta_{\bar{x}})w=a_{ij}w_{\bar{x}^i\bar{x}^i}+b_iw_{\bar{x}^i},
\]
where $a_{ij}$ and $b_i$ are given by~(\ref{eq:ab}).
Furthermore, using~(\ref{eq:simple2}), it is easy to see that 
\[
a_{ij}=\frac{1}{\rho^2}\delta_{ij}-2\frac{\rho_{\bar{x}^i}\bar{x}^j}{\rho^2\Psi}+
\frac{\bar{x}^i\bar{x}^j|\nabla\rho|^2}{\rho^2\Psi^2}.
\]
We now turn to the proof of~(\ref{eq:neumanntheta}). Note that 
\[
u_{x^i}=\nabla_{\bar{x}}w\cdot\bar{x}_{x^i}
=\nabla_{\bar{x}}w\cdot\big(\pi e_i+\pi_{x^i}(x-\a(t))\big)
=\nabla_{\bar{x}}w\cdot\big(\frac{e_i}{\rho}-\frac{\rho_{\bar{x}^i}}{\rho^2\Psi}\rho\bar{x}\big).
\]
From here, we infer the formula
\be\label{eq:grad}
\nabla u(x)=\frac{1}{\rho}\nabla w(\bar{x})
-\frac{\nabla\rho(\bar{x})}{\rho(\bar{x})\Psi(\bar{x})}\bar{x}\cdot\nabla w(\bar{x}).
\ee
From the above formula, we obtain
\be\label{eq:transf}
[\nabla u]^+_-\circ\phi(\xi)=\frac{1}{\rho(\xi)}[\nabla w]^+_-(\xi)
-\frac{\nabla_g\rho(\xi)}{\rho(\xi)\Psi(\xi)}\xi\cdot[\nabla w]^+_-(\xi)
\ee
Note that 
\[
\xi\cdot[\nabla w]^+_-(\xi)=\xi\cdot\big([\nabla_gw]^+_-
+[w_n]^+_-n_{\S}\big)=[w_n]^+_-,
\]
since $[\nabla_gw]^+_-=0$ (due to the fact that $u^{+}|_{\S}=u^{-}|_{\S}$) and
$\xi\cdot n_{\S}=|n_{\S}|^2=1$ ($n_{\S}$ stands for the unit normal on $\S$). 
Observe further that $\rho(\xi)=r(\xi)$. Furthermore,
since $\nabla\rho(\xi)\cdot\xi=\nabla_g\rho\cdot\xi=0$, we
also have $\Psi(\xi)=\rho(\xi)=r(\xi)$.
Form these observations and from~(\ref{eq:transf}) we obtain the formula
\[
[\nabla u]^+_-\circ\phi=\frac{1}{\rho}[w_n]^+_-n_{\S}
-\frac{\nabla_g\rho}{\rho\Psi}[w_n]^+_-
=\frac{1}{r}[w_n]^+_-n_{\S}-\frac{\nabla_gr}{r^2}[w_n]^+_-.
\]
It is straightforward to see that $n_{\Gamma}^i\circ\phi
=\frac{r}{|g|}n^i_{\S}-\frac{\nabla_gf\cdot\nabla_g\xi^i}{|g|}$.
Hence
\[
\begin{array}{l}
\displaystyle
\Big([\nabla u]^+_-\cdot n_{\Gamma}\Big)\circ\phi
=\Big(\frac{1}{r}[w_n]^+_-n_{\S}^i
-\frac{\nabla_gf\cdot\nabla_g\xi^i}{r^2}[w_n]^+_-\Big)
\cdot\Big(\frac{r}{|g|}n_{\S}^i-\frac{\nabla_gf\cdot\nabla_g\xi^i}{|g|}\Big)\\
\displaystyle
=\frac{1}{|g|}[w_n]^+_-
+\frac{|\nabla_gf|^2}{|g|r^2}[w_n]^+_-
=\frac{|g|}{r^2}[w_n]^+_- 
\end{array}
\]
and this proves~(\ref{eq:neumanntheta}).
\prfe

The author wishes to thank Yan Guo for many stimulating discussions. This work has been partly supported by 
the NSF grant CMG 0530862.
%
%

%
\address{ Department of Mathematics,\\
Massachusetts Institute of Technology\\
77 Massachusetts Ave.\\
Cambridge, MA 02139\\
email:{hadzic@math.mit.edu}}


\begin{thebibliography}{}
\bibitem{Lu}
\textsc{Luckhaus, S.}:
Solutions for the two-phase Stefan problem with the Gibbs-Thomson law for the
melting temperature.
{\em Europ. J. Appl. Math.}
{\bf 1}, 101-111 (1990)
\end{thebibliography}
\end{document}